\documentclass[12pt]{article}
\usepackage{amsfonts,amssymb,latexsym}

\newtheorem{theorem}{Theorem}[section]

\newtheorem{proposition}[theorem]{Proposition}

\newenvironment{proof}
{\par\addvspace{0.3cm}\noindent{\rm Proof. }}
{\nopagebreak\mbox{}\hfill $\Box$\par\addvspace{0.25cm}}
\newenvironment{proofof}[1]
{\par\addvspace{0.3cm}\noindent{\rm Proof of #1. }}
{\nopagebreak\mbox{}\hfill $\Box$\par\addvspace{0.25cm}}
\renewcommand{\Re}{\mbox{\rm Re\,}}

\newcommand{\R}{{\mathbb R}}
\newcommand{\C}{{\mathbb C}}
\newcommand{\Z}{{\mathbb Z}}
\newcommand{\T}{{\mathbb T}}
\newcommand{\W}{\mathcal{W}}
\newcommand{\Li}{L^\iy(\T)}

\renewcommand{\kappa}{\varkappa}
\renewcommand{\rho}{\varrho}

\newcommand{\be}{\begin{equation}}
\newcommand{\ee}{\end{equation}}
\newcommand{\ds}{\displaystyle}
\newcommand{\bqn}{\begin{eqnarray}}
\newcommand{\eqn}{\end{eqnarray}}
\newcommand{\nn}{\nonumber}
\newcommand{\ba}{\begin{array}}
\newcommand{\ea}{\end{array}}
\newcommand{\al}{\alpha}
\newcommand{\iv}{^{-1}}
\newcommand{\iy}{\infty}
\newcommand{\twomat}[1]{\left(\ba{cc}#1\ea\right)}
\newcommand{\ta}{\tilde{a}}

\newcommand{\tc}{\tilde{c}}
\newcommand{\bh}{\hat{b}}
\newcommand{\an}{_{\al,n}}

\begin{document}

\date{}
\title{Dyson's constants in the asymptotics of the determinants
of Wiener-Hopf-Hankel operators with the sine kernel}
\author{Torsten Ehrhardt
	\thanks{ehrhardt@math.ucsc.edu.}\\
	Department of Mathematics\\
	University of California\\
	Santa Cruz, CA-95065, USA}
\maketitle

\begin{abstract}
Let  $K_\al^\pm $ stand for the integral operators with the sine kernels
$\frac{\sin(x-y)}{\pi(x-y)}\pm \frac{\sin(x+y)}{\pi(x+y)}$ acting on $L^2[0,\al]$. 
Dyson conjectured that the asymptotics of the 
Fredholm determinants of $I-K_\al^\pm$ are given by
$$
\log\det(I-K_{\al}^\pm) =
-\frac{\alpha^2}{4}\mp \frac{\alpha}{2}-\frac{\log\alpha}{8}+\frac{\log 2}{24}\pm \frac{\log 2}{4}
+\frac{3}{2} \zeta'(-1)+o(1),
\quad\al\to\iy.
$$
In this paper we are going to give a proof of these two asymptotic formulas.
\end{abstract}

\section{Introduction}

In random matrix theory one is interested in the three Fredholm determinants
$$
\det(I-K_\al),\quad \det(I-K_\al^+),\quad 
\det(I-K_\al^-),
$$
where $K_\al$ is the integral operator on $L^2[0,\al]$ with the sine kernel
\bqn
k (x,y)&=&\frac{\sin(x-y)}{\pi(x-y)}
\eqn
and $K_\al^\pm$ are the integral operators on $L^2[0,\al]$ with the Wiener-Hopf-Hankel
sine kernels
\bqn
k^\pm (x,y)&=&\frac{\sin(x-y)}{\pi(x-y)}\pm \frac{\sin(x+y)}{\pi(x+y)}.
\eqn
These determinants are related to the probabilities $E_\beta(n;\alpha)$ that 
in the bulk scaling limit of the three classical Gaussian ensembles of random matrices an interval of lenght $\al$ contains precisely $n$ eigenvalues. It is customary to associate the parameter 
$\beta=2$ with the Gaussian Unitary Ensemble,  $\beta=1$ with the Gaussian Orthogonal
Ensemble,  and $\beta=4$ with the Gaussian Symplectic Ensemble.
The basic relationship between these probabilities and the Fredholm determinants is given by
$$
E_2(0;\alpha)=\det(I-K_\alpha),\quad E_1(0;\alpha)=\det(I-K_\alpha^+),
$$
and
$$
E_4(0;\al)=\frac{1}{2}\left(\det(I-K_{2\al}^+)+\det(I-K_{2 \al}^-)\right),
$$
while expressions for $E_\beta(n;\al)$ with $n\ge1$ also exist \cite{Me,BTW}.

A  problem which has been open for a long time was the rigorous derivation of the asymptotics of these
determinants as $\alpha \to \iy$. Dyson \cite{D1} was able to give a heuristic derivation and conjectured that
\bqn\label{f.Dys}
\log \det(I-K_{2\alpha}) &=&
-\frac{\alpha^2}{2}-\frac{\log\alpha}{4}+\frac{\log 2}{12}+3\zeta'(-1)+o(1),
\qquad \al\to\iy,
\eqn
where $\zeta$ stands for the Riemann zeta function. It is known \cite{Me} that 
\bqn\label{f.4}
\det(I-K_\al^+) = \prod_{n=0}^\iy (1-\lambda_{2n}(\al)),\qquad
\det (I-K_\al^-) = \prod_{n=0}^\iy (1-\lambda_{2n+1}(\al)),
\eqn
where $\lambda_n(\al)$  are the decreasingly ordered eigenvalues of the operator $K_{2\al}$.
Using (\ref{f.Dys}) and a non-rigorous derivation of the asymptotics of 
the quotient 
\bqn\label{D.ev}
\frac{\det(I-K_\al^+)}{\det (I-K_\al^-)} &=& \prod_{n=0}^\iy \frac{1-\lambda_{2n}(\al)}{1-\lambda_{2n+1}(\al)},
\eqn
which was given by des Cloiseaux and Mehta \cite{DM},
Dyson obtained the asymptotics formulas 
\bqn\label{f.Dys.pm}
\log \det(I-K^\pm_\alpha) &=&
-\frac{\alpha^2}{4}\mp\frac{\alpha}{2}-\frac{\log \alpha}{8} +\frac{\log 2}{24}\pm \frac{\log 2}{4}+
\frac{3}{2}\zeta'(-1)+o(1),\quad \al\to\iy.
\eqn

Recently  the asymptotic formula (\ref{f.Dys}) was proved independently by Krasovsky \cite{Kr} and the author \cite{E} using different methods. Yet another proof was given by Deift, Its, Krasovsky, and Zhou \cite{DIKZ}. The proofs \cite{Kr, DIKZ} are based on the Riemann-Hilbert method, while the proof
\cite{E} is based on determinant identities and the asymptotics of Wiener-Hopf-Hankel determinants
with certain Fisher-Hartwig symbols \cite{BE3}.

The goal of this paper is to give a proof of (\ref{f.Dys.pm}). In contrast to Dyson's derivation we will
not rely on (\ref{f.Dys}) and (\ref{D.ev}). In fact, we will use methods similar to those of  \cite{E}. As a consequence of (\ref{f.4}), the asymptotic formulas (\ref{f.Dys.pm}) then imply the asymptotic formula (\ref{f.Dys}). Hence the results of the present paper give a fourth derivation of (\ref{f.Dys}).

As was pointed out to the author by A.~Its, another proof of (\ref{f.Dys.pm}), which is based on the Riemann-Hilbert method, can very likely be accomplished. It would rely on (\ref{f.Dys}) and (\ref{f.4}) and involve a (rigorous) derivation of the asymptotics of (\ref{D.ev}) based on observations made in \cite[p.~205/206]{DIZ}.

Let us conclude this introduction with some remarks on what else is known
about the Fredholm determinants under consideration.
It was shown by Jimbo, Miwa, M\^{o}ri, and Sato
\cite{JMMS} (see also \cite{TW}) that the function
$$
\sigma(\al)=\al\frac{d}{d\al}\log\det (I-K_\al)
$$
satisfies a Painlev\'e V equation.
Widom  \cite{W2,W3}  was able to identify the highest term in the asymptotics
of $\sigma(\al)$ as $\alpha\to\iy$. Knowing these asymptotics
one can derive a complete asymptotic expansion for $\sigma(\al)$. 
By integration it follows that the asymptotics of $\det(I-K_{2\al})$ are given by
\bqn\label{f.6}
\log \det(I-K_{2\alpha}) &=&
-\frac{\alpha^2}{2}-\frac{\log\alpha}{4}+C+
\sum_{n=1}^N \frac{C_{2n}}{\alpha^{2n}}+O(\alpha^{2N+2}),\quad
\al\to\iy,
\eqn
with constants $C_{2n}$ that can be computed recursively. However, the constant $C$ cannot
be obtained in this way, and its rigorous identification was done - as mentioned above - only in \cite{Kr,E,DIKZ}.
The asymptotic formula (\ref{f.6}) was obtained in \cite{DIZ} as well; also, in the earlier work by B. Suleimanov 
\cite{Su} a rigorous derivation of the leading term of the asymptotics of the derivative of $\sigma(\alpha)$ was obtained. 

In a similar way, it turns out that the functions 
$$
\sigma_\pm(\al)=\alpha\frac{d}{d\alpha} \log \det (I-K_\alpha^\pm)
$$
satisfy a Painlev\'e  III equation \cite{TW-b,TW-d}. Moreover, the operators
$K_\alpha^\pm$ are related to special cases of integral operators $K_{\nu,\alpha}$ on $L^2[0,\alpha]$
with Bessel kernel,
$$
k_{\nu}(x,y)=\frac{J_\nu(\sqrt{x})\sqrt{y}J_\nu'(\sqrt{y})-\sqrt{x} J_\nu'(\sqrt{x}) J_\nu(\sqrt{y})}{2(x-y)}, \quad \nu>-1.
$$
In fact, $\det(I-K_{\alpha}^\pm)=\det (I-K_{\mp 1/2,\alpha^2})$. In the Bessel case, functions defined similarly to $\sigma_\pm(x)$ satisfy also a Painlev\'e III equation. The determinants $\det(I-K_{\nu,\alpha})$
are the probabilities that no eigenvalues lie in an interval of length $\alpha$ for the Laguerre or Jacobi
random matrix ensembles in the hard edge scaling limit.

It is also interesting to observe that the following identity between
$\det(I-K_\alpha^\pm)$ and $\det(I-K_\alpha) $ exists (see, e.g., \cite{TW}):
\bqn\label{f.rel}
\log \det(I-K^\pm_\alpha) &=& \frac{1}{2}\log \det (I-K_{2\alpha})\mp\frac{1}{2}
\int_{0}^{\alpha}\sqrt{-\frac{d^2}{dx^2}\log\det(I-K_{2x})}\,dx
\eqn
Using this formula it is possible to derive from (\ref{f.6}) a complete asymptotic
expansion for $\log\det(I-K^\pm_\alpha)$ at infinitiy with the exception of the constant,
which remains undetermined due to the integration. Thus, once (\ref{f.6}) had been proved,
the only open problem was to identify the constant terms in (\ref{f.Dys}) and (\ref{f.Dys.pm}).

Let us shortly outline how the paper is organized.  In the following section we will fix the 
basic notation and make some additional comments about the idea of the proof.
We will follow essentially the same lines as in \cite{E}. 
In fact, the proof is even somewhat simpler since some technical results are not needed here
(namely, Prop.~4.2 and Prop.~4.9 of \cite{E}). The auxiliary results which are needed here are either 
the same as or analogous to those of \cite{E}. 
In Section \ref{s3} we will  prove a formula involving Hankel determinants and in Section \ref{s4} we will finally prove the asymptotic formula (\ref{f.Dys.pm}).

\section{Basic notation and some remarks}

We start with introducing some notation. We will denote the real line by $\R$, the positive real
half-axis  by $\R_+$, and the complex unit circle by $\T$. By $L^p(M)$ we will
denote the Lebesgue spaces ($1\le p\le \iy$), where in our cases $M$ is any of the above sets or
a finite subinterval of $\R$.

The $n\times n$ Toeplitz and Hankel matrices are defined by 
\be
T_n(a)=(a_{j-k})_{j,k=0}^{n-1},\qquad
H_n(a)=(a_{j+k+1})_{j,k=0}^{n-1},
\ee
where $a\in L^1(\T)$ and 
$$
a_k =\frac{1}{2\pi}\int_0^{2\pi}a(e^{i\theta})e^{-ik\theta}\, d\theta,\qquad
k\in\Z,
$$
are its Fourier coefficients.
We will also need differently defined  $n\times n$ Hankel matrices
\be\label{f.Hn2}
H_n[b]=(b_{j+k+1})_{j,k=0}^{n-1},
\ee
where the numbers $b_k$ are the (scaled) moments of a function $b\in L^1[-1,1]$, i.e., 
$$
b_{k}=\frac{1}{\pi}\int_{-1}^1 b(x) (2x)^{k-1}\,dx,\qquad k\ge1.
$$

For $a\in \Li$ the Toeplitz and Hankel operators are bounded linear operators
acting on the Hardy space
$$
H^2(\T)=\Big\{\;f\in L^2(\T)\;:\;f_k=0\mbox{ for all } k<0\;\Big\}
$$
by 
\be\label{f.THdef}
T(a)=PM(a)P|_{H^2(\T)},\qquad H(a)=PM(a)JP|_{H^2(\T)},
\ee
where
$P:\sum_{k=-\iy}^\iy f_k t^k\mapsto\sum_{k=0}^\iy f_k t^k$ stands for the Riesz projection, 
$J:f(t)\mapsto t\iv f(t\iv)$ stands for a flip operator, and
$M(a):f(t)\to a(t)f(t)$ stands for the multiplication operator. (These last three operators 
are acting on $L^2(\T)$.) Finally, introduce the projections
\bqn
P_n: \sum_{k\ge0}f_k t^k \in H^2(\T)\mapsto \sum_{k=0 }^{n-1} f_kt^k\in H^2(\T),
\eqn
the image of which can be naturally identified with $\C^n$. Using this we can make the
identifications $P_n T(a)P_n\cong T_n(a)$, $P_n H(a) P_n\cong H_n(a)$.

We will also need the notion of a trace class operator acting on a Hilbert space $H$.
This is a compact operator $A$ such that the series of its singular $s_n(A)$
(i.e., the eigenvalues of $(A^*A)^{1/2}$ counted according to their algebraic multiplicities)
converges. The class of all trace class operators can be made to a Banach space
by introducing the norm
\bqn
\|A\|_{1} &=& \sum_{n\ge1} s_n(A).
\eqn
This class is also a two-sided ideal in the algebra of all bounded linear operators on $H$.
The importance of trace class operators is that for such operators $A$, the operator
trace ``$\mathrm{trace}(A)$'' and the operator determinant ``$\det (I+A)$'' can be defined as
generalizations of matrix trace and matrix determinant. More detailed information on 
this subject can be found, e.g., in \cite{GK}.

For $a\in L^\iy(\R)$, let $M_\R(a):f(x)\mapsto a(x)f(x)$ stand for the multiplication
operator acting on $L^2(\R)$. The convolution operator $W_0(a)$ 
(or, ``two-sided'' Wiener-Hopf operator) is defined by
$$
W_0(a)=\mathcal{F} M_\R(a)\mathcal{F}\iv,
$$
where $\mathcal{F}$ stands for the Fourier transform on $L^2(\R)$.
The continuous analogues of Toeplitz and Hankel operators are operators defined
\bqn\label{f.WH}
W(a)&=& \Pi_+ W_0(a) \Pi_+|_{L^2(\R_+)},\\[1ex]
\label{f.HaR}
H_\R(a)&=&\Pi_+W_0(a)\hat{J}\Pi_+|_{L^2(\R_+)},
\eqn
where $(\hat{J}f)(x)=f(-x)$, and $\Pi_+=M_\R(\chi_{\R_+})$ is the projection 
operator on the positive real half axis. The operator $W(a)$ is usually called a Wiener-Hopf operator,
and we will refer to $H_\R(a)$ as a Hankel operator, too. (The notation will avoid a possible confusion between $H_\R(a)$ and $H(a)$.)
One can show that if $a\in L^1(\R)$, then
$W(a)$ and $H_\R(a)$ are integral operators on $L^2(\R)$ with the
kernel $\hat{a}(x-y)$ and $\hat{a}(x+y)$, respectively, where
$$
\hat{a}(\xi)=\frac{1}{2\pi}\int_{-\iy}^\iy e^{-ix\xi} a(x)\,dx
$$
stands for the Fourier transform of $a$. For $\al>0$ we will define the projection operator
\bqn\label{f.Pi}
\Pi_\al&:& f(t)\in L^2(\R_+)\mapsto \chi_{[0,\al]}(x)f(x)\in L^2(\R_+).
\eqn
The image of this operator can be identified with $L^2[0,\al]$.

With this notation the integral operators $K_\alpha$ and $K_\alpha^\pm$ can now be
seen to be truncated Wiener-Hopf and Wiener-Hopf-Hankel operators,
$$
K_\alpha=\Pi_\alpha W(\chi)\Pi_\alpha|_{L^2[0,\al]},\qquad
K_\alpha^\pm =\Pi_\alpha( W(\chi)\pm H_\R(\chi))\Pi_\alpha|_{L^2[0,\al]},
$$
where $\chi$ stands for the characteristic function of the interval $[-1,1]$.
Notice that $K_\al$ and $K_\al^\pm$ are trace class operators and that 
\bqn
\det(I-K_\alpha)  &=& \det\Big[\Pi_\alpha W(1-\chi)\Pi_\alpha\Big],\\[1ex]
\det(I-K_\alpha^\pm) &=& \det \Big[\Pi_\alpha\Big(W(1-\chi)\pm H_\R(1-\chi)\Big)\Pi_\alpha\Big].
\eqn
For determinants of Wiener-Hopf operators (and, more recently, also for determinants of  Wiener-Hopf-Hankel 
operators) results describing the asymptotics as $\al\to\iy$ exist under the condition
that the underlying symbol is sufficiently well behaved. These results are known
as  Achiezer-Kac formulas (if the symbol has no singularities) and as Fisher-Hartwig type formulas
(if the symbol has a finite number of certain types of singularities). An overview about this topic can be found in \cite{BS}.
In our case the symbol is the characteristic function $1-\chi$ vanishing on the interval $[-1,1]$,
a state of affairs which is not covered by the just mentioned cases and which renders the situation 
completely non-trivial. 

The main idea of the proof given in this paper is to relate the Fredholm determinants
$\det (I-K_\al^\pm)$  to the determinants of different operators for which Fisher-Hartwig type formulas can be applied. Let us introduce the functions
\bqn\label{f.u-hat}
\hat{u}_\beta(x)=\left(\frac{x-i}{x+i}\right)^\beta,\qquad
\hat{v}_\beta(x)=\left(\frac{x^2}{1+x^2}\right)^\beta,\qquad x\in\R,\;\; \beta\in\C,
\eqn
where these functions are supposed to be continuous on $\R\setminus\{0\}$ and to have
their values approaching $1$ as $x\to\pm \iy$.
Then we are going to prove that 
\bqn\label{21}
\det(I-K_\al^\pm ) &=& \exp\left(-\frac{\al^2}{4}\mp\frac{\al}{2}\right)
\det\Big[\Pi_{\al}(I\pm H_{\R}(\hat{u}_{\mp 1/2}))\iv \Pi_{\al}\Big].
\eqn
It is now illuminating to point out that the determinants on the right hand side can be identified
with determinants of truncated Wiener-Hopf-Hankel operators with Fisher-Hartwig symbols.
In fact, it is proved in \cite{BE3} that 
\bqn
\det\Big[\Pi_\al(W(\hat{v}_\beta)+H_\R(\hat{v}_\beta))\Pi_\al\Big]
&=& e^{-\al\beta} \det \Big[\Pi_\alpha(I+H_\R(\hat{u}_{-\beta}))\iv\Pi_\al\Big]
\quad \mbox{ if } -\frac{1}{2}< \Re \beta <\frac{3}{2},\nn \\
\det\Big[\Pi_\al(W(\hat{v}_\beta)-H_\R(\hat{v}_\beta))\Pi_\al\Big]
&=& e^{-\al\beta} \det \Big[\Pi_\alpha(I-H_\R(\hat{u}_{-\beta}))\iv\Pi_\al\Big]
\quad \mbox{ if } -\frac{1}{2}< \Re \beta <\frac{1}{2}.\nn
\eqn
However, we will avoid making use of these formulas for two reasons. First of all, the determinants
on the right hand side of (\ref{21}) are those occurring primarily in the proof, and their asymptotics are
computed also in \cite{BE3}. Secondly, the left hand side of the last formula is, as it stands, not defined for $\beta=-1/2$. (It can be definined by analytic continuation in $\beta$ because the
right hand side makes sense for $-3/2<\Re\beta<1/2$.)

\section{A Hankel determinant formula}
\label{s3}

In this section we are going to prove two formulas of the kind
\bqn
\det H_n[b] &=& G^n\det\Big[P_n (I+H(\psi))\iv P_n\Big],\nn
\eqn
where $b\in L^1[-1,1]$ is a (sufficiently smooth) continuous and nonvanishing function 
on $[-1,1]$ multiplied in one case with the function $(1+ x)^{1/2}$ and in another case 
with $(1-x)^{-1/2}$.
The function $\psi$ and the constant $G$ depend on $b$. 
A  formula of the same type was proved in \cite{E}. However, 
the conditions on the function $b$ and the form of the function $\psi$ were different.

Before we state the result we have to introduce more notation. 
Let $\W$ stand for the Wiener algebra, i.e., the set of all
functions in $L^1(\T)$ whose Fourier series are absolutely convergent. Moreover, 
let 
\be
\W_\pm= \left\{ \;a\in\W\;:\; a_n=0\mbox{ for all }\pm n <0\;\right\},
\ee
be two Banach subalgebras of $\W$,
where $a_n$ stand for the Fourier coefficients of $a$.
Notice that $a\in\W_+$ if and only if $\ta\in\W_-$, where $\ta(t):=a(t\iv)$, $t\in\T$.
Finally, we denote  by $G\W$ and $G\W_\pm$ the group of invertible elements
in the Banach algebras $\W$ and $\W_\pm$, respectively.

A function $a\in\W$ is said to admit a canonical Wiener-Hopf factorization in $\W$ 
if it can be written in the form
\bqn
a(t)=a_-(t)a_+(t),\qquad t\in\T,
\eqn
where $a_\pm\in G\W_\pm$. It is easy to see that $a\in\W$ 
admits a canonical Wiener-Hopf factorization in $\W$ if and only if $a\in G\W$ and if 
the winding number of $a$ is zero (see, e.g., \cite{BS}) . Moreover, this condition is
equivalent to the existence of a logarithm $\log a$ which belongs to $\W$. If this is
fulfilled, then one can unambiguously define the geometric mean of $a$ by
\bqn
G[a] &:=& \exp\Big(\frac{1}{2\pi}\int_{0}^{2\pi} \log a(e^{i\theta})\, d\theta\Big).
\eqn

The following result (which is not yet what we ultimately want) is cited from \cite[Thm.~4.5]{E}. 
The invertibility statement is taken from \cite[Prop.~4.3]{E}.
Recall that a function $a$ on $\T$ is called even if $\ta=a$, where $\ta(t):=a(t\iv)$.

\begin{theorem}\label{t4.5}
Let $a\in G\W$ be an even function which possesses a canonical Wiener-Hopf factorization
$a(t)=a_-(t)a_+(t)$. Define $\psi(t)=\ta_+(t)a_+\iv(t)$, and let $b\in L^1(\T)$ be
\bqn\label{f.TH-H1}
b(\cos\theta)=a(e^{i\theta})\sqrt{\frac{1+\cos\theta}{1-\cos\theta}}.
\eqn
Then $I+H(\psi)$ is invertible on $H^2(\T)$ and  
\bqn
\det H_n[b] &=&G[a]^n \det\Big[ P_n (I+H(\psi))\iv P_n\Big].
\eqn
\end{theorem}

In order to be able to state the desired result we introduce (for $\tau\in\T$ and $\beta\in\C$)
the functions
\bqn\label{f.u}
u_{\beta,\tau}(e^{i\theta})=\exp(i\beta(\theta-\theta_0-\pi)),
\quad 0<\theta-\theta_0<2\pi,\quad \tau=e^{i\theta_0}.
\eqn
These functions are continuous on $\T\setminus\{\tau\}$ and have a jump discontinuity
at $t=\tau$ whose size is determined by $\beta$.

The promised formulas are now given in the following theorem. Notice that the difference between
Theorem \ref{t4.5} and Theorem \ref{t2.4} (as well as to Thm.~4.6 of \cite{E}) is in the conditions
on the underlying functions.

\begin{theorem}\label{t2.4}
Let $c\in G\W$ be an even function which possesses a canonical Wiener-Hopf factorization
$c(t) = c_-(t)c_+(t)$. Define $b^+,b^-\in L^1[-1,1]$ and $\psi^+,\psi^-\in L^\iy(\T)$ by
$$
\ba{rclcrcl}
b^+(\cos\theta) &= &c(e^{i\theta})\sqrt{2+2\cos\theta} ,   &\qquad&
\psi^+(e^{i\theta}) &=& \tc_+(e^{i\theta}) c_+\iv(e^{i\theta}) u_{-1/2,1}(e^{i\theta}),\\[2ex]
b^-(\cos\theta) &= &\ds\frac{c(e^{i\theta})}{\sqrt{2-2\cos\theta}} ,   &\qquad&
\psi^-(e^{i\theta}) &=& \tc_+(e^{i\theta}) c_+\iv(e^{i\theta}) u_{1/2,-1}(e^{i\theta}).
\ea
$$
Then the operators $I+H(\psi^\pm)$ are invertible on $H^2(\T)$ and
\bqn\label{f.H.H}
\det H_n[b^\pm] &=& G[c]^n\det\Big[P_n (I+H(\psi^\pm))\iv P_n\Big].
\eqn
\end{theorem}

For the proof of this theorem we will apply some auxiliary results, which are stated in \cite{E}
in  connection with Thm.~4.6 and which we are not going to restate here. However, we will
recall the following notation, which is used here and later on.
For $r\in [0,1)$ and $\tau \in \T$ let $G_{r,\tau}$ be the following operator acting on $\Li$:
\be\label{f.G}
G_{r,\tau}:a(t)\mapsto b(t)=a\left(\tau\frac{t+r}{1+rt}\right)
\ee

\begin{proofof}{Theorem \ref{t2.4}} 
The first problem is to verify the invertibility of $I+H(\psi^{\pm})$.
In the special case $c_+\equiv1$, i.e., for $I+H(u_{-1/2,1})$ and $I+H(u_{1/2,-1})$,
this was done in  \cite[Thm.~3.6]{BE3}.
(Notice that $I+H(u_{1/2,-1})$ is similar to $I-H(u_{1/2,1})$.)
The proof in the case where $c_+\not\equiv 1$ can be done in the same way as in \cite[Sec.~3.2]{BE3}
or  \cite[Prop.~4.2]{E}.
We refrain from copying the proof with the little modifications necessary and make instead only the 
following remarks.

The proof in \cite{BE3} consists of two parts. First one determines the essential spectrum of
the Hankel operators. Since $\psi^{\pm}$ have their discontinuities at the same places 
as $u_{\mp1/2,\pm 1}$ and the one-sided limits there are also the same, the essential spectrum of
$H(\psi^\pm)$ is the same as that of $H(u_{\mp 1/2,\pm 1})$.
The second step is to determine the kernel of the operators $I+H(u_{\mp 1/2,\pm 1})$.
(Passing to the adjoints, gives similarly information about the cokernel.) The crucial point in
\cite{BE3} is to write, e.g.,  $u_{\beta,1}=\xi_{-\beta}\eta_{\beta}$ with $\xi_{-\beta}=(1-t\iv)^{-\beta}$, 
$\eta_{\beta}(t)=(1-t)^{\beta}$, $t\in\T$. What one uses about these functions are the facts that 
$\xi_{-\beta}(t)=1/\eta_{\beta}(t\iv)$,  that they and their inverses belong to certain Hardy spaces.
 In our case, one has to write, e.g., 
$\psi^+=(\tc_+ \xi_{1/2})\cdot (c_+\iv \eta_{-1/2})$.
The factors in this product have the just mentioned properties, too. Hence the proof works in the same way.

The proof of (\ref{f.H.H})
will be carried out by an approximation argument and with the help of Theorem \ref{t4.5}.
For $r\in[0,1)$ consider the even functions
$$
a_r^\pm (t)= c(t)\Big((1\mp rt)(1\mp rt\iv)\Big)^{\pm 1/2},\qquad t\in\T.
$$
Clearly, $a_r^\pm \in G\W$.
The functions $b_r^\pm$ defined in terms of $a_r^\pm$ by formula (\ref{f.TH-H1}) 
evaluate to
\bqn
b_r^\pm (x) &=& 
b^\pm (x) (2\pm 2x)^{\mp 1/2}
\Big(1+r^2 \mp 2rx \Big)^{\pm 1/2}
\sqrt{\frac{1+x}{1-x}}\nn\\
&=& 
b^\pm (x)
\left(\frac{2\mp 2 x}{1+r^2 \mp 2rx} \right)^{\mp 1/2},
\qquad\qquad
x\in(-1,1).\nn
\eqn
Then the functions $b_r^\pm$ converge to $b^\pm$ in the norm of $L^1[-1,1]$.
Hence, if we fix $n$,
$$
\det H_n[b^\pm]=\lim_{r\to1}\det H_n[b_r^\pm].
$$

It is now easily seen that the canonical Wiener-Hopf factorization of $a_r^\pm$ is given
by $a_r^\pm(t)=a_{r,-}^\pm(t)a_{r,+}^\pm(t)$
with the factors
$$
a_{r,+}^\pm(t)=c_+(t)(1\mp rt)^{\pm 1/2},\qquad
a_{r,-}^\pm(t)=c_-(t)(1\mp rt\iv)^{\pm 1/2}.
$$
Notice also that $G[a_r]=G[c]$.
If we define
$$
\psi_r^\pm (t):=\ta_{r,+}^\pm (t)(a_{r,+}^\pm(t))\iv=
\tc_{+}(t)c_{+}\iv(t)\left(\frac{1\mp rt}{1\mp rt\iv}\right)^{\mp 1/2},
$$
we  can apply Theorem \ref{t4.5} and conclude that
$$
\det H_n[b_r^\pm]=G[c]^n\det\Big[P_n ( I+H(\psi_r^\pm))\iv P_n\Big].
$$
It follows that 
$$
\det H_n[b^\pm]= G[c]^n\lim_{r\to 1} \det\Big[P_n ( I+H(\psi_r^\pm ))\iv P_n\Big].
$$

Next define 
\bqn
f_r^\pm(t) &:=& 
\left(\frac{1\mp rt}{1\mp rt\iv}\right)^{\mp 1/2}\nn
\eqn
and observe that $f_r^{\pm}\to u_{\mp 1/2,\pm 1}$ in measure as $r\to 1$.
Hence also $\psi_r^\pm \to \psi^\pm$ in measure. 
Because the sequence $\psi_r^\pm$ is bounded in
the $L^\iy$-norm it follows  
that $H(\psi_r^\pm)$ converges strongly to $H(\psi^\pm)$
on $H^2(\T)$ (see, e.g.,  Lemma 4.7 of \cite{E}).
In order to obtain that
\bqn\label{f.sc}
(I+H(\psi_r^\pm))\iv & \to &(I+H(\psi^\pm))\iv
\eqn
strongly on $H^2(\T)$,
it is necessary and sufficient that the following stability condition,
$$
\sup_{r\in[r_0,1)}\left\|(I+H(\psi_r^\pm))\iv\right\| < \iy,
$$
is satisfied (see, e.g., Lemma 4.8 of \cite{E}). Here $r_0$ is some number in $[0,1)$.

Stability criteria for such a type of operator sequences were established in 
\cite{ES1} (see Sections~4.1, 4.2, and 5.2 therein), and we are going to apply the 
corresponding results. First of all, there exist certain mappings
$\Phi_0$ and $\Phi_\tau$, $\tau\in \T$, which are defined by
$$
\Phi_0[\psi_r]:=\mu\mbox{-}\lim\limits_{r\to 1} \psi_r,\qquad
\Phi_0[\psi_r]:=\mu\mbox{-}\lim\limits_{r\to 1} G_{r,\tau}\psi_r.
$$
Here $\mu\mbox{-}\lim$ stands for the limit in measure.
It is now easy to see that these mappings evaluate as follows,
$$
\Phi_0[f_r^\pm ] = u_{\mp 1/2,\pm 1}, \qquad
\Phi_\tau[f_r^\pm] = u_{\mp 1/2,\pm 1}(\tau),
$$
if $\tau\neq\pm 1$, and 
$$
\Phi_{\pm 1}[f_r^\pm]= \mu\mbox{-}\lim_{r\to 1} G_{r,\pm 1} f_r^\pm=
\mu\mbox{-}\lim_{r\to 1}
\left(\frac{1+rt}{1+rt\iv}\right)^{\pm1/2} = u_{\pm 1/2,-1}
$$
if $\tau=\pm 1$. Because of $\psi_r^\pm =\tc_+ c_+\iv f_r^\pm $ it follows immediately that
\bqn
\Phi_0[\psi_r^\pm] &=& \tc_+c_+\iv u_{\mp 1/2,\pm 1},\nn\\[1ex]
\Phi_{\pm 1}[\psi_r^\pm ] &=& u_{\pm 1/2,-1},\nn\\[1ex]
\Phi_\tau[\psi_r^\pm ] &=& \mbox{constant function, }\quad\tau\in \T\setminus\{\pm 1\}.\nn
\eqn
The stability criterion in \cite{ES1} (Thm.~4.2 and Thm.~4.3) says that $I+H(\psi_r^\pm)$ is stable if and only if each of the following operators is invertible:
\begin{itemize}
\item[(i)] $\Psi_0[I+H(\psi_r^\pm)]=I+H(\Phi_0[\psi_r^\pm])=I+H(\psi^\pm)$,
\item[(ii)] $\Psi_{\pm 1}[I+H(\psi_r^\pm)]=I\pm H(\Phi_{\pm 1}[\psi_r^\pm])=I\pm H(u_{\pm 1/2,-1})$,
\item[(iii)] $\Psi_{\mp 1}[I+H(\psi_r^\pm )]=I\mp H(\Phi_{\mp}[\psi_r^\pm ])=I$,
\item[(iv)] $\Psi_\tau[I+H(\psi_r^\pm )]=$
$$
\twomat{I&0\\0&I}+\twomat{P&0\\0&Q}\twomat{M(\Phi_\tau[\psi_r^\pm ])&0\\0&
M(\widetilde{\Phi_{\bar{\tau} }[\psi_r^\pm ]})}\twomat{0&I\\ I&0}
\twomat{P&0\\0&Q}=\twomat{I&0\\0&I}
$$
($\tau\in\T$, $\mathrm{Im}(\tau)>0$)
\end{itemize}
The invertibility is obvious for (iii) and (iv). As to (i) and (ii) the invertibility has been stated at the beginning of the proof. Notice that $I\pm H(U_{\pm 1/2,-1})$ is similar to $I\mp H(u_{\pm 1/2, 1})$.

We can thus conclude that the sequence $I+H(\psi_r^\pm )$ is stable and the strong convergence (\ref{f.sc}) follows. Hence the matrices
$P_n(I+H(\psi_r^\pm ))\iv P_n$ converge to  $P_n (I+H(\psi^\pm))\iv P_n$
as $r\to 1$. This implies that their determinants also converge and proves the assertion.
\end{proofof}

\section{Proof of the asymptotic formula}
\label{s4}

In order to prove the asymptotic formula
(\ref{f.Dys.pm}), we are going to discretize the underlying Wiener-Hopf-plus-Hankel operators 
$I-K_\al^\pm$. This will give us Toeplitz-plus-Hankel operators. Let 
$\chi_{\alpha}$ denote the characteristic function of the subarc
$\{e^{i\theta}:\al<\theta<2\pi-\al\}$ of $\T$.

\begin{proposition}\label{p3.1}
For each $\al>0$ we have
\bqn
\det(I-K_\al^\pm)=\lim_{n\to\iy} 
\det\Big[ T_n(\chi_{\frac{\al}{n}})\pm H_n(\chi_{\frac{\al}{n}})\Big].
\eqn
\end{proposition}
\begin{proof}
The operator $K_\alpha^\pm$ is the integral operator on
$L^2[0,\alpha]$ with the kernel $K(x-y)\pm K(x+y)$, where
$K(x)=\frac{\sin x}{\pi x}$.
Consider the $n\times n$ matrices
\bqn
A_n^\pm &=& \left[\frac{\al}{n}K\left(\frac{\al(j-k)}{n}\right)\pm
\frac{\al}{n}K\left(\frac{\al(j+k+1)}{n}\right)
\right]_{j,k=0}^{n-1},\nn
\\
B_n^\pm &=& \left[\frac{\al}{n}\int_{0}^1\int_0^1
\left\{
K\left(\frac{\al(j-k+\xi-\eta)}{n}\right)\pm
K\left(\frac{\al(j+k+\xi+\eta)}{n}\right)\right\}\,d\xi d\eta
\right]_{j,k=0}^{n-1}.\nn
\eqn
The entries of $A_n^\pm-B_n^\pm$ can be estimated
uniformly by $O(n^{-2})$ using the mean value theorem. Hence the 
 Hilbert-Schmidt norm of $A_n^\pm-B_n^\pm$ is $O(n\iv)$, and the trace norm is $O(1/\sqrt{n})$.

The rest of the proof can be completed in the same way as in \cite[Prop.~5.1]{E} by showing that 
$\det(I-K_\al^\pm )=\det (I_n - B_n^\pm)$ and $\det(I-A_n^\pm)=
\det\Big[ T_n(\chi_{\frac{\al}{n}})\pm H_n(\chi_{\frac{\al}{n}})\Big]$.
\end{proof}

After discretizing, the next goal is to reduce the Toeplitz-plus-Hankel
determinants to Hankel determinants. For this purpose we use an exact identity
which is stated in the following result cited from \cite[Thm.~2.3]{BE1}.

\begin{proposition} \label{p2.3}
Let $a\in L^1(\T)$ be an even function, and let $b\in L^1[-1,1]$ be given by
\bqn\label{f.TH-H}
b(\cos\theta)=a(e^{i\theta})\sqrt{\frac{1+\cos\theta}{1-\cos\theta}}.
\eqn
Then $\det\Big[T_n(a)+H_n(a)\Big]=\det H_n[b]$.
\end{proposition}

Notice that the assumption $b\in L^1[-1,1]$ implies that $a\in L^1(\T)$. Applying 
the previous result yields the following.

\begin{proposition}\label{p3.3}
For each $\al>0$ and $n\in \mathbb{N}$ we have
\bqn\label{f.33}
\det \Big[T_n(\chi_{\frac{\al}{n}})\pm H_n(\chi_{\frac{\al}{n}})\Big] &=& 
(\mu\an)^{\pm n/2} \left(\frac{\rho_{\al,n}+1}{2}\right)^{n^2}\det H_n[b_{\alpha,n}^\pm],
\eqn
where
\be\label{f.bpm}
b_{\alpha,n}^+(x)=
\sqrt{\frac{2+2x}{1+\mu\an^2-2\mu\an x}},\qquad
b_{\alpha,n}^-(x)=
\sqrt{\frac{1+\mu\an^2+2\mu\an x}{2-2x}},\qquad
\ee
and $\rho\an$ and $\mu\an$ are numbers (unambiguosly) defined by
\be\label{f.rho}
\rho\an = \cos\left(\frac{\al}{n}\right),
\qquad
\frac{1+\mu\an^2}{2\mu\an} = \frac{3-\rho_{\al,n}}{1+\rho_{\al,n}},\qquad 0<\mu\an<1.
\ee
\end{proposition}
\begin{proof}
In the plus-case, we apply Proposition \ref{p2.3} with
$$
a(e^{i\theta})=\chi_{\frac{\al}{n}}(e^{i\theta}),\quad 
b(x)=\bh\an^+(x):=\chi_{[-1,\rho_{\alpha,n}]}(x)\sqrt{\frac{1+x}{1-x}}.
$$
In the minus-case, we apply this proposition with 
$$
a(e^{i\theta})=\chi_{\frac{\al}{n}}(-e^{i\theta}),\quad 
b(x)=\bh\an^-(x):=\chi_{[-\rho_{\alpha,n},1]}(x)\sqrt{\frac{1+x}{1-x}}.
$$
Hence we obtain (by using the general formula $\det(T_n(f)+H_n(f))=\det(T_n(\hat{f})-H_n(\hat{f}))$
with $\hat{f}(t)=f(-t)$ in the minus-case)
\bqn
\det \Big[T_n(\chi_{\frac{\al}{n}})\pm H_n(\chi_{\frac{\al}{n}})\Big] &=&
\det H_n[\bh\an^\pm].\nn
\eqn
The entries of $H_n[\bh\an^\pm]$ are the moments $[\bh\an^\pm]_{1+j+k}$, $0\le j,k\le n-1$.
Computing them yields
\bqn
[\bh\an^+]_k  &=& \frac{1}{\pi}\int_{-1}^{\rho_{\al,n}} \sqrt{\frac{1+x}{1-x}} (2x)^{k-1}\,dx
\nn\\
&=&
\frac{1}{\pi}
\left(\frac{\rho_{\al,n}+1}{2}\right)^k\int_{-1}^1
\sqrt{\frac{1+y}{\frac{3-\rho_{\al,n}}{1+\rho_{\al,n}}-y}}
\left(2y-2\,\frac{1-\rho_{\al,n}}{1+\rho_{\al,n}}\right)^{k-1}dy\nn\\
&=&
\frac{\sqrt{\mu\an}}{\pi}
\left(\frac{\rho_{\al,n}+1}{2}\right)^k\int_{-1}^1
b\an^+(y)
\left(2y-2\,\frac{1-\rho_{\al,n}}{1+\rho_{\al,n}}\right)^{k-1}dy\nn
\eqn
and
\bqn
[\bh\an^-]_k  &=& \frac{1}{\pi}\int_{-\rho_{\al,n}}^{1} \sqrt{\frac{1+x}{1-x}} (2x)^{k-1}\,dx
\nn\\
&=&
\frac{1}{\pi}
\left(\frac{\rho_{\al,n}+1}{2}\right)^k\int_{-1}^1
\sqrt{\frac{\frac{3-\rho_{\al,n}}{1+\rho_{\al,n}}+y}{1-y}}
\left(2y+2\,\frac{1-\rho_{\al,n}}{1+\rho_{\al,n}}\right)^{k-1}dy\nn\\
&=&
\frac{1}{\pi\sqrt{\mu\an}}
\left(\frac{\rho_{\al,n}+1}{2}\right)^k\int_{-1}^1
b\an^-(y)
\left(2y+2\,\frac{1-\rho_{\al,n}}{1+\rho_{\al,n}}\right)^{k-1}dy.\nn
\eqn
Hence
$$
H_n[\bh\an^\pm] = \left[(\mu\an)^{\pm 1/2}\left(\frac{\rho\an+1}{2}\right)^{j+k+1}
\frac{1}{\pi} \int_{-1}^1 b\an^\pm(y)(2y\mp 2\tau\an)^{j+k}\,dy\right]_{j,k=0}^{n-1}.
$$
with certain $\tau\an$. One can pull out certain diagonal matrices from the left and the
right, which give the terms $(\mu\an)^{\pm n/2}((1+\rho\an)/2)^{n^2}$ after taking
the determinant. The remaining matrix can be written as
$$
\left[\frac{1}{\pi} \int_{-1}^1 b\an^\pm(y)
(2y\mp 2\tau\an)^j (2y\mp 2\tau\an)^k\,dy\right]_{j,k=0}^{n-1}.
$$
After expanding $(2y\mp 2\tau\an)^j$ and $(2y\mp 2\tau\an)^k$ into two binomial series
it is easily seen that the previous matrix is the matrix $H_n[b\an^\pm]$ multiplied
from the left and right with triangular matrices having ones on the diagonal.
This implies the desired assertion.
\end{proof}

In the following result and also later on we use the functions
\be\label{f.psi.pm1}
\psi_{\al,n}^{\pm } (t):= 
\left(\mp\frac{t\mp\mu\an}{1\mp \mu\an t}\right)^{\mp1/2}
\ee
with the sequence $\mu\an$ defined by (\ref{f.rho}).

\begin{proposition}\label{p3.4}
For each $\al>0$ we have
\bqn\label{f.41}
\lim_{n\to\iy} \det\Big[ T_n(\chi_{\frac{\al}{n}})\pm  T_n(\chi_{\frac{\al}{n}})\Big] &=&
\exp\Big(-\frac{\al^2}{8}\mp \frac{\al}{2}\Big)
\lim_{n\to\iy}
\det\Big[ P_n(I+H(\psi_{\al,n}^\pm))\iv P_n\Big]. 
\qquad
\eqn
\end{proposition}
\begin{proof}
The asymptotics (as $n\to\iy$) of the numbers appearing in (\ref{f.33}) of Proposition \ref{p3.3}
are given by
$$
\frac{1+\rho_{\al,n}}{2}=1-\frac{\al^2}{4 n^2}+O(n^{-4}),\qquad
\mu\an=1-\frac{\al}{n}+O(n^{-2}).
$$
Hence using this proposition it follows that $$
\lim_{n\to\iy} \det \Big[T_n(\chi_{\frac{\al}{n}})\pm H_n(\chi_{\frac{\al}{n}})\Big] =
\exp\Big(-\frac{\al^2}{8}\mp \frac{\al}{2}\Big)
\lim_{n\to\iy}
\det H_n[b_{\al,n}^\pm].
$$
Next we introduce
$$
c(e^{i\theta})=\left((1\mp \mu\an t)(1\mp \mu\an t\iv )\right)^{\mp 1/2},
$$
and we are going to employ Theorem \ref{t2.4} .
It can be verified easily that $G[c]=1$ and that $c(t)=c_-(t)c_+(t)$ is a canonical Wiener-Hopf factorization of $c$ where
$c_+(t)=(1\mp \mu_{\al,n} t)^{\mp 1/2}$ and $c_-(t)=(1\mp \mu_{\al,n} t\iv )^{\mp 1/2}$. 
Moreover,
$$
\tc_+(t)c_+\iv(t)=\left(\frac{1\mp \mu\an t}{1\pm \mu\an t\iv}\right)^{\pm 1/2}.
$$
The functions  $b^\pm$ and $\psi^\pm$ defined in Theorem \ref{t2.4} now evaluate to
\bqn
b^\pm(x)  &=& (1+\mu\an^2 \mp 2\mu\an x)^{\mp1/2}(2\pm 2x)^{\pm 1/2}
\;=\;b\an^\pm(x),\nn\\[1ex]
\psi^\pm(t) &=& \left(\frac{1\mp\mu\an t}{1\mp\mu\an t\iv}\right)^{\pm 1/2}(\mp t)^{\mp1/2}
\;=\;\psi^\pm\an (t).\nn
\eqn
Combining all this we obtain from Theorem \ref{t2.4} that 
$$
\det H_n[b_{\al,n}^\pm ] =\det\Big[P_n(I+H(\psi_{\al,n}^\pm ))\iv P_n\Big],
$$
which concludes the proof.
\end{proof}

The next step is to identify the limit on the right hand side of (\ref{f.41}). For this purpose
we resort to an auxiliary result, which was stated in \cite{E} (with a slight change of notation).
In order to make the reference correct, we allow (for the time being) $\mu_{\al,n}\in [0,1)$ 
to be an arbitrary sequence and define the functions
\be\label{f.hal}
h_{\alpha}(t) = \exp\left(-\al\frac{1-t}{1+t}\right),
\qquad
h_{\al,n}(t)  = \left(\frac{t+\mu_{\al,n}}{1+\mu_{\al,n}t}\right)^n.
\ee
Moreover, we also consider the functions
$\psi\an^\pm$ as being defined by (\ref{f.psi.pm1}) with this arbitrary sequence.

\begin{proposition}\label{p3.5}
Let $\alpha>0$ be fixed, and consider (\ref{f.psi.pm1}) and (\ref{f.hal}). Assume that 
\bqn\label{f.mu.an1}
\mu_{\al,n}&=&1-\frac{\al}{n}+O(n^{-2})\quad
\mbox{ as } n\to\iy.
\eqn
Then the following is true:
\begin{itemize}
\item[(i)]
The operators $H(\psi_{\al,n}^\pm)$ are
unitarily equivalent to the operators
$\pm H(u_{\mp 1/2,1})$.
\item[(ii)]
The operator
$$
P_n(I+H(\psi_{\al,n}^\pm))\iv P_n-P_n
$$
is unitarily equivalent to operators
$$
A_n = H(h_{\alpha,n})(I\pm H(u_{\mp 1/2,1}))\iv H(h_{\alpha,n})-H(h_{\alpha,n})^2,
$$
which are trace class operators and converge as $n\to\iy$ in the trace norm to 
$$
A=H(h_{\alpha})(I\pm H(u_{\mp 1/2,1}))\iv H(h_{\alpha})-H(h_{\alpha})^2.
$$
\end{itemize}
\end{proposition}
\begin{proof}
These results are proved in \cite[Prop.~4.12]{E} with a change of the condition on the sequence
$\mu\an$. This change is consistent with the different notation for $h_\al$.
In fact, one has only to replace $\alpha$ by $2\alpha$.
Moreover, instead of the functions
$\psi_\alpha^\pm$, functions $\psi_\alpha^\pm-1$ occur, which do not change the
Hankel operators. The fact that the operators $I\pm H(u_{\mp 1/2,1})$ are invertible 
has already been stated in the proof of Theorem \ref{t2.4} (see also \cite[Prop.~4.1]{E} or \cite[Thm.~3.6]{BE3}).
\end{proof}

In the following proposition we identify the limit of the determinant
appearing in the right hand side of (\ref{f.41}). 
We return to the specific definition of $\mu\an$ given in (\ref{f.rho}) 
and to the definitions (\ref{f.psi.pm1}) and (\ref{f.hal}) in terms of this sequence.

\begin{proposition}\label{p3.5x}
For each $\al>0$ we have
\bqn
\lim_{n\to\iy}
\det \Big[P_n(I+H(\psi_{\al,n}^\pm) )\iv P_n\Big]
&=&
\det \Big[H(h_\al)(I \pm H(u_{\mp 1/2,1}))\iv H(h_\al)\Big].
\label{f.36}
\eqn
\end{proposition}
\begin{proof}
Proceeding as in \cite[Prop.~5.5]{E} we notice that $H(h_\al)^2$ is a projection operator.
(The slight change in notation, $\al\mapsto 2\al$, does not affect the statements made here).
Hence, in the same way it is established that
$$
\det \Big[H(h_\al)(I\pm H(u_{\mp 1/2,1}))\iv H(h_\al)\Big] =
\det \Big[I+H(h_\al)(I\pm H(u_{\mp 1/2,1}))\iv H(h_\al)-H(h_\al)^2\Big].
$$
where the determinant on the left hand side is of that an operator acting on the image of 
$H(h_\al)^2$, while the right hand side corresponds to an operator acting on $L^2(\R_+)$.

Similarily, the determinant on the left hand side of (\ref{f.36}) can be written as
\bqn
\det \Big[P_n(I+H(\psi_{\al,n}^\pm))\iv P_n\Big] &=&
\det \Big[I+P_n(I+H(\psi_{\al,n}^\pm))\iv P_n-P_n\Big]\nn\\[1ex]
&=& \det\Big[I+H(h\an)(I\pm H(u_{\mp 1/2,1}))\iv H(h\an)-H(h\an)^2\Big].\nn
\eqn
As stated in the proof of Proposition \ref{p3.4}, the sequence $\mu\an$ has the asymptotics
(\ref{f.mu.an1}).  By applying the previous proposition the desired assertion follows.
\end{proof}

We are now finally able to identify the determinants $\det(I-K_\al^\pm )$. Recall 
in this connection the definition (\ref{f.u-hat}) of the functions
$\hat{u}_\beta$.

\begin{theorem}
For each $\al>0$ we have
\bqn
\det(I-K_\al^\pm ) &=& \exp\left(-\frac{\al^2}{4}\mp\frac{\al}{2}\right)
\det\Big[\Pi_{\al}(I\pm H_{\R}(\hat{u}_{\mp 1/2}))\iv \Pi_{\al}\Big].
\label{f.64}
\eqn
\end{theorem}
\begin{proof}
{}From Propositions \ref{p3.1}, \ref{p3.4} and \ref{p3.5x} it follows that
\bqn
\det(I-K_\al^\pm) 
&=& \exp\left(-\frac{\al^2}{4}\mp\frac{\al}{2}\right)
\det \Big[H(h_\al)(I\pm H(u_{\mp 1/2,1}))\iv H(h_\al)\Big].\nn
\eqn
As noted in the proof of \cite[Thm.~5.6]{E}, there exists a unitary transform $S$ from
$H^2(\T)$ onto $L^2(\R_+)$ such that $H_{\R}(a)= SH(b)S^* $ with $a(x)=b\left(\frac{1+ix}{1-ix}\right)$.
In the specific case we obtain
$$
H_\R(\hat{u}_{\pm 1/2})=SH(u_{\pm 1/2,1})S^*,\qquad
H_\R( e^{ix\alpha})=S H(h_\al) S^*.
$$
This together with the remark that $H( e^{ix\alpha})^2=\Pi_{\al}$
implies (\ref{f.64}).
\end{proof}

Finally we are using result of \cite{BE3} in order to establish 
the promised asymptotic formula. Recall that $\zeta$ stands for the
Riemann zeta function.

\begin{theorem}
The following asymptotic formula holds as $\al\to\iy$:
\bqn\label{f.as1}
\log\det(I-K_{\al}^\pm) &=& -\frac{\al^2}{4}\mp\frac{\al}{2}
-\frac{\log\al}{8}+\frac{\log 2}{24}\pm \frac{\log 2}{4}+\frac{3}{2}\zeta'(-1) +o(1)
\eqn
\end{theorem}
\begin{proof}
In Sect.~3.6 of \cite{BE3} it has been proved that 
\bqn
\det\Big[\Pi_{\al}(I\pm H_{\R}(\hat{u}_{\mp1/2}))\iv \Pi_{\al}\Big]
&\sim&
\al^{-1/8}\pi^{1/4} 2^{\pm 1/4}G(1/2),\qquad \al\to\iy,
\eqn
where  $G(z)$ stands for the Barnes $G$-function \cite{Bar}.
Notice that $G(3/2)=G(1/2)\Gamma(1/2)$, $\Gamma(1/2)=\pi^{1/2}$,
and $G(1)=1$.
This together with the previous theorem implies that 
\bqn
\log \det(I-K_{\al}^\pm) = -\frac{\al^2}{2}\mp\frac{\al}{2}-\frac{\log\al}{8}
+\frac{\log \pi}{4} \pm \frac{\log 2}{4} +\log G(1/2)+o(1).
\eqn
Finally observe that
\bqn
\log G(1/2) = -\frac{\log\pi}{4}+\frac{3}{2}\zeta'(-1)+\frac{\log 2}{24},\nn
\eqn
which follows from a formula for $G(1/2)$ in terms of Glaisher's constant
$A= \exp(-\zeta'(-1)+1/12)$ given in \cite[page 290]{Bar}.
\end{proof}

\end{document}